 \documentclass[12pt]{article}

 \usepackage{mathrsfs,amsfonts}
 \usepackage[dvips]{color}

 \newtheorem{theorem}{Theorem}[section]
 \newtheorem{lemma}[theorem]{Lemma}
 \newtheorem{corol}[theorem]{Corollary}
 \newtheorem{prop}[theorem]{Proposition}
 \newtheorem{example}[theorem]{Example}
 
 \newtheorem{con1}[theorem]{Condition}
 \newtheorem{remark}[theorem]{Remark}

 \def\blemma{\begin{lemma}\sl{}\def\elemma{\end{lemma}}}
 \def\bproposition{\begin{prop}\sl{}\def\eproposition{\end{prop}}}
 \def\btheorem{\begin{theorem}\sl{}\def\etheorem{\end{theorem}}}

 \def\beqlb{\begin{eqnarray}}\def\eeqlb{\end{eqnarray}}
 \def\beqnn{\begin{eqnarray*}}\def\eeqnn{\end{eqnarray*}}

 \def\proof{\noindent{\it Proof.~~}}\def\qed{\hfill$\Box$\medskip}

 \def\<{\langle}\def\>{\rangle}

 \def\mcr{\mathscr}\def\mbb{\mathbb}

 \def\ar{\!\!&}

\begin{document}

\

\bigskip\bigskip

\centerline{\Large\bf Nonlinear branching processes with immigration }
\

\bigskip

\centerline{Pei-Sen Li}

\centerline{School of Mathematical Sciences, Beijing Normal University,}

\centerline{Beijing 100875}

\centerline{E-mail: \tt peisenli@mail.bnu.edu.cn}

\bigskip

The nonlinear branching process with immigration is constructed as the pathwise unique solution of a stochastic integral equation driven by Poisson random measures. Some criteria for the regularity, recurrence, ergodicity and strong ergodicity of the process are then established.     




\section{Introduction}
Markov branching processes are models for the evolution of populations of particles. Those processes constitute one of the most important subclasses of continuous-time Markov chains. Standard references on those processes are \cite{Harris} and \cite{Athreya72}. The basic property of an ordinary linear branching process is that different particles act independently when giving birth or death. In most realistic situations, however, this property is unlikely to be appropriate. In particular, when the number of particles becomes large or the particles move with high speed, the particles may interact and, as a result, the birth and death rates can either increase or decrease. Those considerations have motivated the study of nonlinear branching processes.
On the other hand, a branching process describes a population evolving randomly in an isolated environment. A useful and realistic modification of the model is the addition of new particles from outside sources. This consideration has provided the stimulation for the study of branching models with immigration and/or resurrection.

Let $\{r_i: i\ge 0\}$ be a sequence of nonnegative constants with $r_0=0$ and $\{b_i: i\ge 0\}$ a discrete probability distribution {on $\mbb{N}:=\{0,1,\ldots\}$ }with $b_1=0$. A continuous-time Markov chain is called a nonlinear branching process if it has density matrix $R=(r_{ij})$ given by
 \beqlb\label{1.4}
r_{ij}=\left\{
\begin{array}{lcl}
r_i b_{j-i+1} & & {j\geq i+1, i\ge 1,}\cr
-r_i & & {j=i\ge 1,}\cr
r_i b_0 & & {j=i-1, i\ge 1,}\cr
0 & & \mbox{otherwise.}
\end{array} \right.
 \eeqlb
A typical special case is where $r_i=\alpha i^\theta$ for $\alpha\ge 0$ and $\theta>0$, which reduces to the ordinary linear branching process when $r_i=\alpha i$.
Let $\gamma\ge 0$ and let $\{a_i: i\ge 0\}$ be another discrete probability distribution on $\mbb{N}$ satisfying $a_0=0$. A continuous-time Markov chain is called a nonlinear branching process with resurrection if its density matrix is given by
 \beqlb\label{1.1}
\rho_{ij}=\left\{
\begin{array}{lcl}
r_i b_{j-i+1} & & {j\geq i+1, i\ge 1,}\cr
-r_i & & {j=i\ge 1,}\cr
r_i b_0 & & {j=i-1, i\ge 1,}\cr
\gamma a_j & & {j> i=0,}\cr
-\gamma & & {j = i=0,}\cr
0 & & \mbox{otherwise.}
\end{array} \right.
 \eeqlb
Here the resurrection means that at each time when the process gets extinct, some immigrants come into the population at rate $\gamma$ according to the distribution $\{a_i\}$. By a nonlinear branching process with immigration we mean a Markov chain with density matrix $Q=(q_{ij})$ given by
 \beqlb\label{1.2}
q_{ij}=\left\{
 \begin{array}{lcl}
r_i b_{j-i+1}+\gamma a_{j-i} & & {j\geq i+1, i\ge 0,}\cr
-r_i-\gamma & & {j=i\ge 0,}\cr
r_i b_0 & & {j=i-1, i\ge 1,}\cr
0 & & \mbox{otherwise.}
 \end{array} \right.
 \eeqlb
In this model, the immigrants come at rate $\gamma$ according to the distribution $\{a_i\}$ independently of the inner population.

The purpose of this paper is to investigate the construction and basic properties of the nonlinear branching process with immigration defined by (\ref{1.2}). Let
 \beqnn
m = \sum_{j=0}^\infty ja_j, \quad M = \sum_{j=0}^\infty jb_j,
 \eeqnn
which represent the birth mean and immigration mean of the process, respectively.
Moreover, we introduce the
functions
 \beqnn
F(s)=\sum_{i=0}^\infty a_is^i, ~ A(s)=\gamma(1-F(s)), ~ G(s)=\sum_{i=0}^\infty b_is^i, ~B(s)= G(s)-s, \qquad s\in [0,1].
 \eeqnn
Let $q$ be the smaller root of the equation $G(s)=s$ in $[0,1]$. We sometimes denote $r_i$ by $r(i)$ for notational convenience.

Suppose that $(\Omega,\mcr{F},\mcr{F}_t,P)$ is a probability space satisfying the usual hypotheses. { Denote $m(i)=b_i$ and $n(i)=a_i$ for each $i\in \mbb{N}$.}
Let $\{p(t)\}$ and $\{q(t)\}$ be $(\mcr{F}_t)$-Poisson point processes with characteristic measures $dum(dz)$ and $\gamma n(dz)$, respectively. We assume $\{p(t)\}$ and $\{q(t)\}$ are independent of each other. Let $N_p(ds,du,dz)$ and $N_q(ds,dz)$ be the Poisson random measures associated with $\{p(t)\}$ and $\{q(t)\}$, respectively.
Given an $\mbb{N}$-valued $\mcr{F}_0$-measurable random variable $X_0,$ let us consider the stochastic integral equation
 \beqlb\label{e0.1}
X_t= X_0+\int_0^t \int_0^{r({X_{s-}})}\int_\mbb{N}(z-1) N_p(ds,du,dz)+  \int_0^t \int_\mbb{N}  z N_q(ds,dz).
 \eeqlb
Let $\zeta=\lim_{k\rightarrow\infty} \tau_k,$ where $\tau_k=\inf\{t\geq0: X_t\geq k\}.$ The above equation only makes sense for $0\leq t< \zeta.$ We call $\zeta$ the explosion time of $\{X_t\}$ and make the convention $X_t=\infty$ for $t\geq \zeta.$ We say the solution is non-explosive if $\zeta= \infty.$ As a special case of (\ref{e0.1}) we also consider the equation
  \beqlb\label{0.2}
X_t= X_0+\int_0^t \int_0^{r({X_{s-}})}\int_\mbb{N}(z-1) N_p(ds,du,dz).
 \eeqlb

We now state the main results of the paper.

\btheorem\label{t0.2a}
There exists a pathwise unique solution to {\rm(\ref{e0.1})}. Moreover, if the solution to {\rm(\ref{0.2})} is non-explosive, then so is the solution to {\rm(\ref{e0.1})}.
 \etheorem

\btheorem\label{t0.4} Let $\{X_t\}$ be the solution to {\rm(\ref{e0.1})} and let $Q_{ij}(t) = P(X_t=j|X_0=i).$ Then  $Q_{ij}(t)$ solves the Kolmogorov forward equation of $Q$.
\etheorem

\btheorem\label{t0.2aa}
The solution to {\rm(\ref{e0.1})} is the minimal process of $Q$ and the solution to {\rm(\ref{0.2})} is the minimal process of $R$.
\etheorem

\btheorem\label{t0.5}
The density matrix $R$ is regular if and only if $Q$ is regular.
\etheorem

\btheorem\label{t1.4}
{\rm(1)} If $M\leq 1,$ then $Q$ is regular.

{\rm(2)} Suppose that $\sum^\infty_{i=1}r_i^{-1}<\infty$. Then $Q$ is regular if and only if $M\leq 1.$

{\rm(3)} Suppose that $1<M\leq \infty$ and $r_i=\alpha i^{\theta}$ for $\alpha>0$ and $\theta>0$. Then $Q$ is regular if and only if for some $\varepsilon\in (q,1)$,
we have
 $$
\int^1_\varepsilon\frac{1}{B(s)}\bigg(\ln \frac{1}{s}\bigg)^{\theta-1}ds=-\infty.
 $$
\etheorem

In the following three theorems, we assume $\gamma r_i b_0>0$ for every $i\ge 1$, so the matrix $Q$ is irreducible.

\btheorem\label{t3.1}
{\rm(1)} Suppose that $m< \infty$, $M< 1$ and $\lim_{i\rightarrow\infty}r_i=\infty$. Then the nonlinear branching process with immigration is recurrence.

{\rm(2)} Suppose that $r_i$ is increasing and there exist constants $\alpha>0$ and $N>0$ such that $r_i/i\geq \alpha$ holds for each $i> N$. Then the nonlinear branching process with immigration is recurrent if $M\leq 1$ and
 \beqnn
J := \int_0^1 \frac{1}{\alpha B(y)}\cdot \exp\bigg[-\int^y_0 \frac{A(x)}{\alpha B(x)}dx\bigg]dy=\infty.
 \eeqnn

{\rm(3)} Suppose that $M>1$. Then the nonlinear branching process with immigration is transient.

{\rm(4)} Suppose that $r_i$ is increasing and there exist constants $\alpha>0$ and $N>0$ such that $r_i/i\leq \alpha$ holds for each $i> N$. Then the nonlinear branching process with immigration is transient if $M\leq 1$ and
 \beqnn
J := \int_0^1 \frac{1}{\alpha B(y)}\cdot \exp\bigg[-\int^y_0 \frac{A(x)}{\alpha B(x)}dx\bigg]dy<\infty.
 \eeqnn
 \etheorem

\btheorem\label{t4.1}
{\rm(1)} If $m< \infty$, $M\leq 1$, $r_i$ is increasing and $\sum_{i=1}^\infty r_i^{-1}< \infty$, then the nonlinear branching process with immigration is ergodic.

{\rm(2)} Suppose that $r_i=\alpha i^\theta$ for $\alpha>0$ and $\theta\ge 1$. Then the recurrent nonlinear branching process with immigration is ergodic if and only if
 \beqlb\label{4.8}
\int_0^1 \frac{A(s)}{\alpha B(s)}\bigg(\ln \frac{1}{s}\bigg)^{\theta-1} ds < \infty.
 \eeqlb
{\rm(3)} If $m< \infty, M<1$ and $\lim\inf_{i\rightarrow\infty} r_i/i>0$, then the nonlinear branching process with immigration is exponentially ergodic.
\etheorem

\btheorem\label{t5.1} {\rm(1)} If $m< \infty$, $M< 1$, $r_i$ is increasing and $\sum_{i=1}^\infty r_i^{-1}< \infty$, then the process is strongly ergodic.

{\rm(2)} Suppose that $r_i=\alpha i^\theta$ for $\alpha>0$ and $\theta> 1$. Then the nonlinear branching process with immigration is strongly ergodic if and only if
 \beqlb\label{5.1}
\int^1_0 \frac{1}{\alpha B(s)}\bigg(\ln\frac{1}{s}\bigg)^{\theta-1}ds < \infty.
 \eeqlb

{\rm(3)} If $\sum_{i=1}^\infty r_i^{-1}= \infty$, then the nonlinear branching process with immigration is not strongly ergodic.
\etheorem

The nonlinear branching process with resurrection defined above was introduced by \cite{ChenRR97}, who studied the problems of uniqueness, recurrence and ergodicity of the process. The model has attracted the attention of a number of authors. In particular, \cite{Zhang01} gave criteria for strong ergodicity of the process. \cite{Chen05} and \cite{Pakes} established some criteria for their regularity and uniqueness. \cite{Chen02} studied some interesting differential-integral equations associated with a special class of nonlinear branching processes and gave some characterizations of their mean extinction times. \cite{Chen06} established a Harris regularity criterion for such processes. The existence and uniqueness of linear branching processes with instantaneous resurrection were studied in \cite{Chen90}. However, most of the study of models with immigration have been focused on linear branching structures. The branching process with immigration was studied in \cite{Karlin}, who gave a characterization of the one-dimensional marginal distributions of the process starting from zero. An ergodicity criterion for the process was given in \cite{Yang}. \cite{LJP06} established some recurrence criteria for linear branching processes with immigration and resurrection.

The first three theorems above give constructions of nonlinear branching processes with and without immigration.
These provide convenient formulations of the processes. In particular, the result of Theorem~\ref{t0.5} is derived as an immediate consequence of (\ref{e0.1}) and (\ref{0.2}). We hope the equations can also be useful in some other similar situations. The proof of Theorem~\ref{t1.4} is based on Theorem~\ref{t0.5} and the results of \cite{ChenRR97} and \cite{Chen06}.

The study of recurrence of the immigration model is more delicate since the problem cannot be reduced to the extinction problem of the original nonlinear branching process as in the case of a resurrection model. Theorem~\ref{t3.1} was proved by using the results of the minimal nonnegative solutions as developed in \cite{ChenMF04} and  comparing the process with some linear branching processes which was studied by \cite{LJP06}.

The proofs of the ergodicities in Theorems~\ref{t4.1} and~\ref{t5.1} are based on comparisons of the process with some suitably designed birth-death process and estimates of the mean extinction time.

\section{Stochastic integral equations}

\setcounter{equation}{0}

Stochastic integral equations with jumps have been playing increasingly important roles in the study of Markov processes. In this section, we give a construction of the solution to (\ref{e0.1}) and prove the solution is a minimal nonlinear branching process with immigration. This result is then used to study the regularity of the density matrix $Q$. We refer to \cite{Ikeda89} for the general theory of stochastic equations with jumps.

 \bproposition\label{t0.0}  The pathwise uniqueness of solutions holds for the equation {\rm(\ref{e0.1})}.
 \eproposition
 \proof Let $\{X_t\}$ and $\{X'_t\}$ be any two solutions of equation (\ref{e0.1}) with $X_0=X'_0$.
 By passing to the conditional probability $P(\cdot|\mcr{F}_0),$ we may and do assume $X_0= X'_0$ is deterministic.
 Let $\tau_m=\inf\Big\{t\geq 0: X_t \geq m\Big\},$ $\tau'_m=\inf\Big\{t\geq 0: X'_t \geq m\Big\}$ and $\sigma_m=\tau_m\wedge \tau'_m.$ It is sufficient to show that $\tau_m=\tau'_m=\sigma_m$ and $X_t=X'_t$ for all $t\leq \sigma_m~(m=1,2,\ldots).$ Then
\beqnn
X_{t\wedge\sigma_m} - X'_{t\wedge\sigma_m}
 \ar=\ar
 \int^{t\wedge\sigma_m}_0\int^\infty_0\int_{\mbb{N}_{m+1}}
 (z-1)[1_{\{0 < u\leq r(X_{s-})\}}\cr
 \ar\ar\qquad\qquad\qquad\qquad
-\, 1_{\{0<u\leq r(X'({s-}))\}}]N_p(ds,du,dz),
 \eeqnn
 where $\mbb{N}_m=\{0,1,2,\ldots,m\}.$ Taking the expectation, we get
 \beqnn
\ar\ar E[|X_{t\wedge\sigma_m} - X'_{t\wedge\sigma_m}|]\cr
 \ar\ar\qquad
\leq E\bigg\{\int^{t\wedge\sigma_m}_0\int^\infty_0\int_{\mbb{N}_{m+1}}
|(z-1)[1_{\{0 < u\leq r(X_{s-})\}} - 1_{\{0 < u\leq r(X'_{s-})\}}]|N_p(ds,du,dz)\bigg\}\cr
 \ar\ar\qquad
\leq E\bigg\{\int^{t\wedge\sigma_m}_0\int^\infty_0\int_{\mbb{N}_{m+1}}
(z+1)|1_{\{0 < u\leq r(X_{s-})\}} - 1_{\{0 < u\leq r(X'_{s-})\}}|dsdum(dz)\bigg\}\cr
 \ar\ar\qquad
\leq (M_{m+1}+1)E\bigg\{\int^{t\wedge\sigma_m}_0
|r(X_{s-})-r(X'_{s-})|ds\bigg\}\cr
 \ar\ar\qquad
\leq (M_{m+1}+1)\int^{t}_0
E[|r(X_{s\wedge\sigma_m-})-r(X'_{s\wedge\sigma_m-})|]ds,
\eeqnn
where $M_m:=\int_{\mbb{N}_m}z m(dz).$ By taking $m\geq X_0,$ we have $X_{s-}\vee X'_{s-} \leq m$ for $0<s\leq \sigma_m.$
 Denote $d_m=\sup\{|(r(i)-r(j))/(i-j)|:i\neq j,\quad 0\leq i,j\leq m\}.$ Then we have
\beqlb\label{0.3}
\ar\ar E[|X_{t\wedge\sigma_m} - X'_{t\wedge\sigma_m}|]\cr
 \ar\ar\qquad
\leq
(M_{m+1}+1)d_m\int^{t}_0
E[|X_{s\wedge\sigma_m-}-X'_{s\wedge\sigma_m-}|]ds.
\eeqlb
Since $X_{s\wedge\sigma_m}$ and $X'_{s\wedge\sigma_m}$ only have countably many discontinuous points, we can also use $X_{s\wedge\sigma_m}$ and $X'_{s\wedge\sigma_m}$ instead of $X_{s\wedge\sigma_m-}$ and $X'_{s\wedge\sigma_m-}$ in the right hand side of (\ref{0.3}).
Using Gronwall's inequality we have $E[|X_{s\wedge\sigma_m}-X'_{s\wedge\sigma_m}|]=0$. Thus we can conclude that $X_t=X'_t$ for all $t\in [0,\sigma_m)$ a.s. This clearly implies that $\tau_m=\tau'_m= \sigma_m$ a.s. and the pathwise uniqueness of solutions of (\ref{e0.1}) is proven.

\btheorem\label{t0.2} For any $\mbb{N}$-valued $\mcr{F}_0$-measurable random variable $X_0,$ there is a pathwise unique solution to~{\rm(\ref{0.2})}.
\etheorem
\proof Without loss of generality, we assume $X_0$ is deterministic. Let $D_1=\{s: p(s)\in (0,r(X_0)]\times \mathbb{N}\}.$ Since
$$
E[N_p((0,t]\times(0,r(X_0)]\times\mbb{N})]=
\int^t_0ds\int^{r(X_0)}_0du\int_{\mathbb{N}}m(dz)=
tr(X_0)<\infty,
$$  the set $D_1$ is discrete in $(0,\infty)$. Let $\sigma_1$ be the minimal  element in $D_1$ and $p(\sigma_1)=(u_1,z_1)$.
Then set
 \beqnn
X_t= \left\{
\begin{array}{lcl}
X_0, & & t\in[0,\sigma_1)\cr
X_0+ (z_1-1), & & t=\sigma_1.
\end{array} \right.
 \eeqnn
The process $\{X_t: 0< t\leq \sigma_1\}$ is clearly the solution of (\ref{0.2}). Set $D_2=\{s: p(s+\sigma_1)\in [0,r(X(\sigma_1))]\times \mathbb{N}\},$   $\sigma_2$ be the minimal  element in $D_2$ and $p(\sigma_1+\sigma_2)=(u_2,z_2)$. Define $\{X_t: \sigma_1< t\leq \sigma_1+\sigma_2\}$ by
\beqnn
X_t= \left\{
\begin{array}{lcl}
x(\sigma_1), & & t\in(\sigma_1,\sigma_1+\sigma_2)\cr
x(\sigma_1)+ (z_2-1), & & t=\sigma_1+\sigma_2.
\end{array} \right.
 \eeqnn
It is easy to see that $\{X_t: 0< t\leq \sigma_1+\sigma_2\}$ is the unique solution of (\ref{0.2}).
Continuing this process successively, we get a process $\{X_t: 0\leq t< \tau\}$, where $\tau=\sum^\infty_{i=1}\sigma_i.$ Next, we show $\tau=\zeta:=\lim_{k\rightarrow\infty}\tau_k,$ where $\tau_k=\inf\{t\geq 0: X_t\geq k\}.$ Clearly, for each $n\geq 0$ we have $X_t< \infty$ for  $t\in [0, \sum^n_{i=0} \sigma_i].$ Then $\sum^n_{i=0} \sigma_i< \zeta$ holds for each $n\geq 0,$ and so $\tau\leq \zeta.$
On the other hand, since
$$
E\bigg[\int^{t\wedge\tau_m}_0\int^{r(X_{s-})}_0\int_{\mbb{N}}N_p(ds, du, dz)\bigg]\leq t\max_{0\leq k\leq m}r(k)<\infty,
$$
the process $\{X_t\}$ has finitely many jumps before $t\wedge\tau_m$, therefore $t\wedge\tau_m< \tau,$ since $t\geq 0$ and $m\geq 1$ can be arbitrary, we get $\zeta\leq \tau.$ Then we have $\tau=\zeta.$
Hence $X_t$ is determined in the time interval $[0, \zeta);$ the uniqueness is clear from Proposition~\ref{t0.0}.
 \qed

\noindent\emph{Proof of  Theorem~{\rm\ref{t0.2a}}.}~
Let $\{X^0_t\}$ denote the solution to (\ref{0.2}). Let $\{v_k: k=1,2,\ldots\}$ be the set of jump times of the Poisson process
 $$
 t\longmapsto \int^t_0\int_{\mbb{N}} N_q(ds,dz).
 $$
 We have clearly $v_k\longrightarrow\infty$ as $k\longrightarrow\infty.$ For $0\leq t<v_1$ set $X_t=X^0_t.$ Suppose that $X_t$ has been defined for $0\leq t< v_k$ and let
 $$
\xi= X_{v_k-}+\int_{\{v_k\}}\int_{\mbb{N}} z N_q(ds,dz).
 $$
Here and in the sequel we make the convention $\infty+\cdots= \infty.$ By the assumption there is also a solution $\{X_t^k\}$ to
 $$
X_t=\xi+\int^t_0\int^{r(X_{s-})}_0\int_{\mbb{N}} (z-1)N_p(v_k+ds,du,dz).
$$
Let $\eta_k$ be the explosion time of $\{X_t^k\}$. If $v_k+\eta_k> v_{k+1}$, we define
$X_t=X_{t-v_k}^k$ for $v_k\leq t<v_{k+1}.$  If $v_k+\eta_k\leq v_{k+1}$,  we set $X_t=X_{t-v_k}^k$ for $v_k\leq t<v_k+\eta_k$ and $X_t=\infty$ for $v_k+\eta_k\leq t< v_{k+1}$. By induction that defines a process $\{X_t\}$, which is clearly the pathwise unique solution to (\ref{e0.1}). Obviously, if the solution of (\ref{0.2}) is non-explosive for each deterministic initial state $X_0=i\in \mbb{N}$, we have $\eta_k= \infty$
for all $k\in\mbb{N}$, and so $\{X_t\}$ is non-explosive. \qed

\noindent\emph{Proof of  Theorem~{\rm\ref{t0.4}}.}~ Let $\tilde{N}_p(ds, du, dz) = N_p(ds, du, dz)-dsdum(dz)$ and $\tilde{N}_q(ds, dz)=\tilde{N}_q(ds, du, dz)-dsdun(dz).$ For any bounded function $f$ on $\mbb{N}$ we have,
 \beqlb\label{0.5}
f(X_{t\wedge\tau_m})
 \ar=\ar
f(X_0)+ \int_0^{t\wedge\tau_m} \int_0^{r(X_{s-})}\int_\mbb{N}[f(X_{s-}+z-1)-f(X_{s-})] N_p(ds,du,dz)\cr
 \ar\ar+  \int_0^{t\wedge\tau_m} \int_\mbb{N} [f(X_{s-}+z)-f(X_{s-})] N_q(ds, dz)\cr
 \ar=\ar
f(X_0)+ \int_0^{t\wedge\tau_m} \int_0^{r(X_{s-})}\int_\mbb{N}[f(X_{s-}+z-1)-f(X_{s-})] ds du m(dz)\cr
 \ar\ar+  \int_0^{t\wedge\tau_m} \int_\mbb{N} [f(X_{s-}+z)-f(X_{s-})]  \gamma ds n (dz) +  M_t(f),
 \eeqlb
where
 \beqnn
M_t(f) \ar:=\ar \int_0^{t\wedge\tau_m} \int_0^{r(X_{s-})}\int_\mbb{N}[f(X_{s-}+z-1)-f(X_{s-})]\tilde{N}_p(ds,du,dz) \cr
 \ar\ar\qquad
+ \int_0^{t\wedge\tau_m} \int_\mbb{N} [f(X_{s-}+z)-f(X_{s-})] \tilde{N}_q(ds, dz)
 \eeqnn
is a martingale. Since $X_s\neq X_{s-}$ for at most countably many $s\geq 0$, we can also use $X_s$ instead of $X_{s-}$ in the right hand side of~(\ref{0.5}). In particular, for $f=1_{\{j\}}$ we have
 \beqnn
1_{\{X_{t\wedge\tau_m}=j\}}\ar=\ar 1_{\{X_0=j\}}+\sum_{k=0}^\infty b_k \int_0^{t\wedge\tau_m} r(X_s)[1_{\{X_s+k-1=j\}}-1_{\{X_s=j\}}] ds\cr
 \ar\ar
+\sum_{k=1}^\infty   \gamma a_k \int_0^{t\wedge\tau_m} [1_{\{X_s+k=j\} }-1_{\{X_s=j\}}]ds +M_t(1_{\{j\}}).
 \eeqnn
Write $E_i = E(\cdot|X_0=i)$ for $i\in \mbb{N}$. Taking the expectation in both sides of the above equation and letting $m\longrightarrow \infty$ we get
 \beqnn
E_i(1_{\{X_{t\wedge\zeta}=j\}})\ar=\ar E_i(1_{\{X_0=j\}}) + \sum_{k=0}^\infty b_k E_i\bigg(\int_0^{t\wedge\zeta} r(X_s)[1_{\{X_s+k-1=j\}}-1_{\{X_s=j\}}] ds\bigg)\cr
 \ar\ar
+ \sum_{k=1}^\infty \gamma a_k E_i\bigg(\int_0^{t\wedge\zeta} [1_{\{X_s+k=j\} }-1_{\{X_s=j\}}]ds\bigg).
 \eeqnn
Obviously, here we can remove the truncation ``$\wedge\zeta$'' and obtain
 \beqnn
Q_{ij}(t)\ar=\ar \delta_{ij}+\sum_{k=0}^{j} b_k \int^t_0 [r_{j-k+1} Q_{i,j-k+1}(s)-r_jQ_{ij}(s)]ds \cr
 \ar\ar
+\sum_{k=1}^j \int_0^t   \gamma a_k [Q_{i,j-k}(s)-Q_{ij}(s)]ds\cr
 \ar=\ar \delta_{ij}+ \int^t_0 \bigg(\sum_{k=1}^{j+1} Q_{ik}(s)r_{k}b_{j-k+1} -Q_{ij}(s)r_j\bigg)ds\cr
 \ar\ar
+\int_0^t \bigg( \sum_{k=0}^{j-1}  Q_{ik}(s)  \gamma a_{j-k}-  \gamma Q_{ij}(s)\bigg)ds.
 \eeqnn
Differentiating both sides we get
 \beqnn
Q_{ij}'(t)
 \ar=\ar
\sum_{k=1}^{j+1} Q_{ik}(t)r_{k}b_{j-k+1} - Q_{ij}(t)r_j\cr
 \ar\ar\qquad
+\sum_{k=0}^{j-1} Q_{ik}(t) \gamma  a_{j-k}-  \gamma Q_{ij}(t)\cr
 \ar=\ar
 \sum_{k=0}^{\infty} Q_{ik}(t)q_{kj}.
 \eeqnn
This is just the Kolmogorov forward equation of $Q$. \qed

\noindent\emph{Proof of  Theorem~{\rm\ref{t0.2aa}}.}~
By Theorem~\ref{t0.2a}, the solution $\{X_t\}$ to (\ref{e0.1}) is a time homogeneous Markov process with state space $\bar{\mbb{N}} := \{0,1,2,\dots,\infty\}$. Suppose that $\sigma_1$ and $z_1$ are given in the proof of Theorem \ref{t0.2}. Let $q(v_1)=y_1$. By the properties of Poisson point process, we can see that $P(\sigma_1>t)=e^{-r(X_0)t},$ $P(z_1=i)=m(\{i\})=b_i,$ $P(v_1>t)=e^{-\gamma t},$ $P(y_1=i)= n(\{i\})= a_i  $  and $\sigma_1$, $z_1$, $v_1,$ $y_1$ are mutually independent. Write $P_i(\cdot) = P(\cdot|X_0=i)$ for $i\in \mbb{N}$. Let $\xi_t=\max\{n+m: \sum^n_{i=0}\sigma_i, \sum^m_{i=0} v_i \leq t \}$. Obviously we have $P_i[X_t=j, \xi_t=0]=\delta_{ij}.$ By the Markov property of $\{X_t\},$
 \beqnn
\ar\ar P_i\bigg\{X_t=j, \xi_t=m+1\bigg\}
 =
P_i\bigg\{1_{\{\sigma_1\wedge v_1<t\}}P_{X_{\sigma_1\wedge v_1}}\bigg[X_{t-\sigma_1\wedge v_1}=j, \xi_{t-\sigma_1\wedge v_1}=m\bigg]\bigg\}\cr
 \ar\ar\qquad
=
P_i\bigg\{1_{\{\sigma_1< t\}}1_{\{v_1\geq \sigma_1\}}P_{X_{\sigma_1}}\bigg[X_{t-\sigma_1}=j, \xi_{t-\sigma_1}=m\bigg]\bigg\}\cr
 \ar\ar\qquad\qquad
+\,P_i\bigg\{1_{\{v_1<t\}}1_{\{v_1< \sigma_1\}}P_{X_{v_1}}\bigg[X_{t-v_1}=j, \xi_{t-v_1}=m\bigg]\bigg\}\cr
 \ar\ar\qquad
=
P_i\bigg\{\int^t_0r_ie^{-r_i(t-s)}e^{-\gamma (t-s)}P_{X_{t-s}}[X_s=j, \xi_s=m]ds\bigg\}\cr
 \ar\ar\qquad\qquad
+\,P_i\bigg\{\int^t_0\gamma e^{-\gamma (t-s)}e^{-r_i (t-s)}P_{X_{t-s}}[X_s=j, \xi_s=m]ds\bigg\}\cr
 \ar\ar\qquad
=P_i\bigg\{\int^t_0r_ie^{-(r_i+\gamma)(t-s)}\sum^\infty_{k=i-1}P(z_1=k-i+1)P_k[X_s=j, \xi_s=m]ds\bigg\}\cr
 \ar\ar\qquad\qquad
+\,P_i\bigg\{\int^t_0\gamma e^{-(r_i+\gamma)(t-s)}\sum^\infty_{k=i+1}P(y_1=k-i)P_k[X_s=j, \xi_s=m]ds\bigg\}\cr
 \ar\ar\qquad
=
\sum_{k\neq i}\int^t_0e^{-(r_i+\gamma)(t-s)}q_{ik}P_k[X_s=j, \xi_s=m]ds.
 \eeqnn
Notice that $$P_i[X_t=j]=\sum^\infty_{m=0}P_i[X_t=j, \xi_t=m].$$
From the theory of Markov chains we know $P_{ij}(t):=P_i[X_t=j]$ is the minimal solution to the Kolmogorov equation of the density matrix $Q$, see Chen (2004, p.78). Then $\{X_t\}$ is the minimal process of the density matrix $Q.$
\qed

\noindent\emph{Proof of Theorem~{\rm\ref{t0.5}}.}~
Suppose that $R$ is regular. Then the minimal solution of its Kolmogorov backward equation is honest i.e. the minimal process of $R$ is non-explosive. Applying Theorems \ref{t0.2a} and  \ref{t0.2aa} we know the minimal process of $Q$ is non-explosive. Thus $Q$ is regular. Conversely, suppose that $R$ is not regular.
Then by Theorem 2.7\,(3) in \cite{Anderson91} there exists a non-trivial solution $(u^*_i) $ to
 \beqnn
u_i\leq\sum_{k\neq i}\frac{ r_{ik}}{2\gamma+r_i}u_k,\qquad 0\leq u_i\leq 1.
 \eeqnn
Since $r_{ik}\leq q_{ik},$ we see $(u^*_i)$ is also a solution to
 $$u_i\leq\sum_{k\neq i}\frac{ q_{ik}}{\gamma+q_i}u_k.$$
Using Theorem 2.7\,(3) in \cite{Anderson91} again, we see $Q$ is not regular.\qed

\noindent\emph{Proof of Theorem~{\rm\ref{t1.4}}.}~  By Theorem~\ref{t0.5} we derive the results from Theorem~1.2 of \cite{ChenRR97} and Theorem~2.3 of \cite{Chen06}. \qed

\section{Recurrence}

\setcounter{equation}{0}

\noindent\emph{Proof of Theorem~{\rm\ref{t3.1}}.}~
(1) Under the assumption, there exists a constant $N\ge 1$ such that $r_i\geq \gamma m/(1-M)$ holds for each $i\ge N$.
Take $x_i=i$ for $i\geq 0$. For $i\geq N$ we have
 \beqnn
\sum^{\infty}_{j=0} q_{ij}x_j
 \ar=\ar
r_ib_0(i-1) + \sum^{\infty}_{j=1} (r_ib_{j+1}+\gamma a_j)(i+j) \cr
 \ar=\ar
(r_i+\gamma)i + r_i(M-1)+\gamma m\leq (r_i+\gamma)i = -q_{ii}x_i.
 \eeqnn
Let $(\pi_{ij})$ be the embedded chain of $(q_{ij})$. The above calculations imply that $(x_i)$ is a finite solution of
 \beqnn
\sum^{\infty}_{j=0} \pi_{ij}x_j\le x_i, \qquad i\geq N.
 \eeqnn
Then $Q$ is recurrent by Theorem 4.24 in \cite{ChenMF04}.

(2) Suppose that $M\le 1$ and $J=\infty$. We shall prove the process is recurrent by comparison arguments. Let $\bar{Q} = (\bar{q}_{ij})$ be the density matrix defined
by \beqnn
\bar{q}_{ij}= \left\{
\begin{array}{lcl}
 \alpha ib_{j-i+1} + \gamma a_{j-i} & & {j \geq i+1}\\
- \alpha i - \gamma  & & {j= i}\\
 \alpha ib_0 & & {j= i-1}\\
0 & & \mbox{otherwise},
\end{array} \right.
 \eeqnn
which corresponds to a linear branching process with
immigration. It was proved in \cite{LJP06} that this process is recurrent.
Next, we define the density matrix $Q^*= (q_{ij}^*)$ by
 \beqnn
q_{ij}^*=\left\{
\begin{array}{lcl}
r_i b_{j-i+1} + \gamma a_{j-i}r_i/\alpha i & & {j \geq i+1}\\
-r_i - \gamma r_i/\alpha i & & {j= i}\\
r_ib_0 & & {j= i-1}\\
q_{ij}& & {i< N}\\
0 & & \mbox{otherwise}.
\end{array} \right.
 \eeqnn

Let $(\bar{\pi}_{ij})$ and $(\pi_{ij}^*)$ denote the embedded chains of $(\bar{q}_{ij})$ and $(q_{ij}^*)$, respectively. It is easy to see that $\bar{\pi}_{ij}= \pi_{ij}^*$ for $i\ge N$ and $j\ge 0$. Then $Q^*$ is also recurrent. For $l\geq i> N$ we have
 \beqnn
\sum_{j=i}^\infty q_{ij} = -r_i b_0\leq \sum_{j=i}^\infty q_{lj}^*.
 \eeqnn
Moreover, we have
 \beqnn
\sum_{j=k}^\infty q_{ij} = \sum_{j=k}^\infty q_{lj}^*=0, \qquad k\leq i-1
 \eeqnn
and
 \beqnn
\sum_{j=k}^\infty q_{ij}\leq\sum_{j=k}^\infty q_{ij}^*
 \leq
\sum_{j=k}^\infty q_{lj}^*, \qquad k\geq l+1.
 \eeqnn
Then $Q$ and $Q^*$ are stochastically comparable, so we can construct a $Q$-process $(X_t)$ and a $Q^*$-process $(X^*_t)$ on some probability space in such a way that $X_0=X^*_0$ and $X_t\le X^*_t$ for all $t\ge 0$; see Example 5.51 in \cite{ChenMF04}. Now the recurrence of $(X_t)$ follows from that of $(X^*_t)$.

(3) Since $M>1$, there exists a $s\in (0,1)$ such that $B(s)<0, \quad i.e. \quad \sum^\infty_{i=0} b_i s^{i-1}<1$. Take $H=\{0\}$ and $x_i=1-s^i.$
For $i\geq 1$ we have
 \beqnn
\sum^\infty_{k=0} \pi_{ik}x_k\ar=\ar \pi_{i, i-1}x_{i-1}+\sum^\infty_{k=1}\pi_{i, i+k}x_{i+k}\cr
\ar=\ar
\frac{r_ib_0}{r_i+\gamma}x_{i-1}+\sum^\infty_{k=1}\frac{r_i b_{k+1}+\gamma a_k}{r_i+\gamma}x_{i+k}\cr
\ar=\ar
\frac{1}{r_i+\gamma}\bigg[r_ib_0(1-s^{i-1})+\sum^\infty_{k=1}\gamma a_k(1-s^{i+k})+\sum^\infty_{k=1}r_ib_{k+1}(1-s^{i+k})\bigg]\cr
\ar=\ar
1-\frac{s^i}{r_i+\gamma}\bigg[r_i\sum^\infty_{k=0}  b_k s^{k-1} + \gamma\sum^\infty_{k=1} a_k s^k\bigg]\geq 1-s^i=x_i.
 \eeqnn
Then the process is transient by Theorem~8.0.2 in \cite{Meyn09}.

(4) Since the proof is similar to that of (2), we omit it.
 \qed

\section{Mean extinction time}

\setcounter{equation}{0}

In this section, we assume $r_i=\alpha i^\theta$ for $\alpha>0$ and $\theta\ge 1$.
Let $(X_t)$ be a realization of the nonlinear branching process with
immigration. Its jump times are given successively by $\tau_0=0$ and $\tau_n=\inf\{t:t>\tau_{n-1}, X_t\neq X_{\tau_{n-1}}\}.$
We also define $\sigma_k=\inf\{t\geq \tau_1: X_t=k\}$. In order to prove the criterion for the ergodicity of $(X_t)$, let us consider the absorbing process $\tilde{X}_t := X_{t\wedge \sigma_0}$. The density matrix of this process is given by:
 \beqnn
\tilde{q}_{ij}=\left\{
 \begin{array}{lcl}
q_{ij} & & {i\neq 0}\\
0 & & {i= 0}.
 \end{array} \right.
 \eeqnn
For this process, we define $\tilde{\tau}_0=0$, $\tilde{\tau}_n=\inf\{t:t>\tilde{\tau}_{n-1}, \tilde{X}(t)\neq \tilde{X}(\tilde{\tau}_{n-1})\}$ and $\tilde{\sigma}_k=\inf\{t\geq \tau_1: \tilde{X}_t=k\}$.
It is easy to see that
 \beqlb\label{4.2}E_i \sigma_0=E_i \tilde{\sigma}_0.
 \eeqlb
Let $(\tilde{p}_{ij}(t))$ and $(\tilde{\phi}_{ij}(\lambda))$ denote the transition function and the resolvent of $(\tilde{X}_t)$, respectively.

\blemma\label{l4.3} For any $i\geq 0$ and $s\in [0,1)$, we have
 \beqlb\label{4.4}
\sum_{j=0}^\infty \tilde{p}_{ij}'(t)s^j = \alpha B(s)\sum_{j=1}^\infty \tilde{p}_{ij}(t)j^\theta s^{j-1} - A(s)\sum_{j=1}^\infty \tilde{p}_{ij}(t)s^j, \qquad t\ge 0,
 \eeqlb
and
 \beqlb\label{4.10}
\lambda\sum_{j=0}^\infty \tilde{\phi}_{ij}(\lambda)s^j-s^i = \alpha B(s)\sum_{j=1}^\infty \tilde{\phi}_{ij}(\lambda) j^\theta s^{j-1} - A(s)\sum_{j=1}^\infty \tilde{\phi}_{ij}(\lambda)s^j, \quad \lambda>0.
 \eeqlb
\elemma

\proof From the Kolmogorov forward equation of the transition function we obtain that
 \beqnn
\tilde{p}_{ij}'(t)=\sum_{k=1}^{j-1} \tilde{p}_{ik}(t)(r_k b_{j-k+1}+\gamma a_{j-k}) - \tilde{p}_{ij}(t)(r_j+\gamma) +  \tilde{p}_{i,j+1}(t) r_{j+1} b_0.
 \eeqnn
Multiplying $s^j$ on both sides of the above equality and then summing over $j$, we have
 \beqnn
\sum_{j=0}^\infty \tilde{p}_{ij}'(t)s^j
 \ar=\ar
\sum_{j=0}^\infty\sum_{k=1}^{j-1}\tilde{p}_{ik}(t)r_k b_{j-k+1}s^j + \gamma\sum_{j=0}^\infty \sum_{k=1}^{j-1}\tilde{p}_{ik}(t) a_{j-k} s^j\cr
 \ar\ar
+\, \sum_{j=0}^\infty \tilde{p}_{i, j+1}(t)r_{j+1} s^j b_0 - \sum_{j=1}^\infty \tilde{p}_{ij}(t) r_j s^j - \gamma \sum_{j=0}^\infty \tilde{p}_{ij}(t) s^j,
 \eeqnn
Then we can interchange the order of summation to see
$$\sum_{j=0}^\infty\sum_{k=0}^{j-1}\tilde{p}_{ik}(t)r_k b_{j-k+1}s^j=\sum_{k\neq l} \tilde{p}_{ik}(t)r_k s^{k-1} \sum_{j=k+1}^\infty b_{j-k+1} s^{j-k+1}$$and
$$\gamma\sum_{j=0}^\infty \sum_{k=0}^{j-1}\tilde{p}_{ik}(t) a_{j-k} s^j= \gamma \sum_{k\neq l} \tilde{p}_{ik}(t)s^k\sum_{j=k+1}^\infty a_{j-k}s^{j-k}.$$
It follows that
 $$\sum_{j=0}^\infty \tilde{p}_{ij}'(t)s^j= \sum_{j=1}^\infty \tilde{p}_{ij}(t)r_j s^{j-1}\alpha B(s)-\sum_{j=0}^\infty
\tilde{p}_{ij}(t)s^j A(s).$$
That proves (\ref{4.4}) and (\ref{4.10}) is just the Laplace transform of (\ref{4.4}).\qed

\blemma\label{l4.4} For any $i, k\geq 1,$ we have $\int_0^\infty \tilde{p}_{ik}(t)dt <\infty$
and $\lim_{t\to\infty} \tilde{p}_{ik}(t)= 0$. Furthermore, for $i\geq 1$ and $s\in [0,1)$, we have
 \beqlb\label{4.5}
\sum_{k=1}^\infty\bigg(\int_0^\infty \tilde{p}_{ik}(t)dt\bigg) s^k < \infty.
 \eeqlb
\elemma

\proof Fixing an $i\geq 1$, we can use the Kolmogorov forward equation to see
 \beqnn
\tilde{p}_{i0}(t)=b_0\alpha\int_0^t \tilde{p}_{i1}(u)du,
 \eeqnn
which means that
 \beqnn
\int_0^\infty \tilde{p}_{i1}(t)dt\leq b_0^{-1}\alpha^{-1} <\infty.
 \eeqnn
Suppose that $\int_0^\infty \tilde{p}_{ik}(t)dt <\infty$ for $k\leq j$. By the Kolmogorov forward equations we can see for $j\geq 1$,
 \beqnn
\tilde{p}_{ij}(t)-\delta_{ij}
 \ar=\ar
\sum_{k=1}^{j-1} (\alpha k^\theta b_{j-k+1}+\gamma a_{j-k})\int_0^t \tilde{p}_{ik}(u)du - (\alpha j^\theta+\gamma) \int_0^t \tilde{p}_{ij}(u)du\cr
 \ar\ar\qquad
+\, \alpha(j+1)^\theta b_0 \int_0^t \tilde{p}_{ij+1}(u)du.
 \eeqnn
Letting $t\to\infty$, we have
 \beqnn
\int_0^\infty \tilde{p}_{ij+1}(t)dt <\infty.
 \eeqnn
Then $\int_0^\infty \tilde{p}_{ik}(t)dt <\infty$ by induction. Since the limit $\lim_{t\to\infty} \tilde{p}_{ik}(t)$ always exists, we see $\lim_{t\to\infty} \tilde{p}_{ik}(t)= 0$ immediately.

We next tend to prove (\ref{4.5}). Since $M\leq 1$, we have $B(s)>0$ for a fixed $s\in[0,1)$. Then there exists a $k\ge 1$ so that $k\alpha B(s)-sA(s)> 0$. Using (\ref{4.4}), we have
 \beqlb\label{4.9}
\sum_{j=0}^\infty \tilde{p}'_{ij}(u)s^j \ar=\ar \alpha B(s)\sum_{j=1}^\infty \tilde{p}_{ij}(u)j^\theta s^{j-1}- A(s)\sum_{j=1}^\infty \tilde{p}_{ij}(u)s^j\cr
 \ar\geq\ar
\alpha B(s)\sum_{j=k+1}^\infty \tilde{p}_{ij}(u)j^\theta s^{j-1} - A(s)\sum_{j=1}^\infty \tilde{p}_{ij}(u)s^j\cr
 \ar\geq\ar
[k \alpha B(s)-sA(s)] \sum_{j=k+1}^\infty \tilde{p}_{ij}(u) s^{j-1} - A(s) \sum_{j=1}^k \tilde{p}_{ij}(u)s^j.
 \eeqlb
Let $\|A\| = \max_{s\in[0,1]}|A(s)|$ and $\|B\| = \max_{s\in[0,1]}|\alpha B(s)|$. Then for each $s\in[0,1)$,
 \beqnn
\int_0^t \sum_{j=0}^\infty |\tilde{p}'_{ij}(u)s^j|du
 \ar\leq\ar
\|B\| \int_0^t\sum_{j=1}^\infty \tilde{p}_{ij}(u)j^\theta s^{j-1}du+\|A\|\int_0^t \sum_{j=1}^\infty \tilde{p}_{ij}(u) s^j du\cr
 \ar\leq\ar
t\|B\|\sum_{j=1}^\infty j^\theta s^{j-1} + t\|A\|\sum_{j=1}^\infty s^j<\infty.
 \eeqnn
Then we use Fubini's theorem to see
 \beqnn
\int_0^t\sum_{j=0}^\infty \tilde{p}'_{ij}(u)s^j du=\sum_{j=0}^\infty \int_0^t \tilde{p}'_{ij}(u)s^jdu.
 \eeqnn
Integrating both sides of (\ref{4.9}),
 \beqnn
\sum_{j=0}^\infty \tilde{p}_{ij}(t)s^j-s^i
 \ar\geq\ar
[k \alpha B(s)-sA(s)]\cdot \sum_{j=k+1}^\infty\bigg(\int_0^t \tilde{p}_{ij}(u)du\bigg)s^{j-1}\cr
 \ar\ar
- A(s)\cdot \sum_{j=1}^k \bigg(\int_0^t \tilde{p}_{ij}(u)du\bigg)s^j .
 \eeqnn
Letting $t\to\infty$ and using the fact that $\int_0^\infty \tilde{p}_{ik}(t)dt <\infty,$ we have
 \beqnn
\sum_{j=k+1}^\infty \bigg(\int_0^\infty \tilde{p}_{ij}(u)du\bigg)s^{j-1}<\infty,
 \eeqnn
which implies (\ref{4.5}).\qed

\bproposition\label{l4.5}
Suppose that the nonlinear branching process with immigration is recurrent and {\rm(\ref{4.8})} holds. Then for $i\ge 1$ we have
 \beqlb\label{4.15}
E_i(\sigma_0)\leq \frac{1}{\Gamma(\theta)} \int^1_0 \frac{1-y^i}{\alpha B(y)}\bigg(\ln\frac{1}{y}\bigg)^{\theta-1}dy\cdot \exp\bigg[\frac{1}{\Gamma(\theta)}\int_0^1 \frac{A(y)}{\alpha B(y)} \bigg(\ln\frac{1}{y}\bigg)^{\theta-1} dy\bigg].
 \eeqlb
and
 \beqlb\label{4.16}
E_i(\sigma_0)\geq \int^1_0 \frac{1-y^i}{\alpha B(y)}\bigg(\ln\frac{1}{y}\bigg)^{\theta-1}dy.
 \eeqlb
\eproposition

\proof
Multiplying (\ref{4.10}) by $(\ln(s/y))^{\theta-1}$, dividing by $\alpha B(s)$ and integrating both sides we have
 $$
\int^s_0 \sum^\infty_{j=1} \tilde{\phi}_{ij}j^\theta y^{j-1} (\ln \frac{s}{y})^{\theta-1}dy=\int_0^s \frac{(\lambda+A(y))\sum_{j=1}^\infty \tilde{\phi}_{ij}(\lambda)y^j-y^i+ \lambda\tilde{\phi}_{i0}(\lambda)}{\alpha B(y)}\bigg(\ln\frac{s}{y}\bigg)^{\theta-1}dy.
 $$
 Letting $y=se^{-\frac{x}{j}}$ in the left hand side of the above equation we get
 $$
\int^s_0 \sum^\infty_{j=1} \tilde{\phi}_{ij}j^\theta y^{j-1} (\ln \frac{s}{y})^{\theta-1}dy
=
\int^\infty_0 \sum_{j=1}^\infty \tilde{\phi}_{ij}(\lambda) s^j x^{\theta-1} e^{-x} dx= \Gamma(\theta) \sum^\infty_{j=1} \tilde{\phi}_{ij}(\lambda) s^j.
 $$
 Using the above two equations we obtain
  \beqlb\label{4.12}
\sum_{j=1}^\infty \tilde{\phi}_{ij}(\lambda)s^j= \frac{1}{\Gamma(\theta)}\int_0^s \frac{(\lambda+A(y))\sum_{j=1}^\infty \tilde{\phi}_{ij}(\lambda)y^j-y^i+ \lambda\tilde{\phi}_{i0}(\lambda)}{\alpha B(y)}\bigg(\ln\frac{s}{y}\bigg)^{\theta-1}dy.
 \eeqlb
For $i\geq 1$, $\lambda>0$ and $s\in[0,1]$ let
 \beqnn
\psi_i(\lambda,s)=\sum_{j=1}^\infty \tilde{\phi}_{ij}(\lambda)s^j.
 \eeqnn
Note that
 \beqnn
\lambda \tilde{\phi}_{i0}(\lambda)= \int^\infty_0 e^{-t} p_{i0}\bigg(\frac{t}{\lambda}\bigg)dt \leq \int^\infty_0 e^{-t}dt = 1.
 \eeqnn
Then, by (\ref{4.12}),
 \beqlb\label{4.17}
\psi_i(\lambda,s)
 \ar\leq\ar
\frac{1}{\Gamma(\theta)} \int_0^1 \frac{1-y^i}{\alpha B(y)}\bigg(\ln\frac{1}{y}\bigg)^{\theta-1}dy\cr
 \ar\ar
+\, \frac{1}{\Gamma(\theta)} \int_0^s \frac{(\lambda+A(y))\psi_i(\lambda,y)}{\alpha B(y)}\bigg(\ln \frac{1}{y}\bigg)^{\theta-1}dy.
 \eeqlb
By Lemma~\ref{l4.4},
 \beqnn
\lim_{\lambda\to0}\lambda\sum_{j=1}^\infty \tilde{\phi}_{ij}(\lambda) =\lim_{\lambda\to0} \sum_{j=1}^\infty \int^\infty_0e^{-t}\tilde{p}_{ij}(\frac{t}{\lambda})dt =0.
 \eeqnn
It follows that, for $s\in[0, 1]$,
 \beqlb\label{4.6}
\lim_{\lambda\to0} \lambda \psi_i(\lambda,s)\leq\lim_{\lambda\to0}\lambda\sum_{j=1}^\infty \tilde{\phi}_{ij}(\lambda) =0.
 \eeqlb
 Denote
 \beqnn
C_i \ar:=\ar \frac{1}{\Gamma(\theta)} \int_0^1 \frac{1-y^i}{\alpha B(y)}\bigg(\ln\frac{1}{y}\bigg)^{\theta-1}dy\\
\ar\leq\ar
\frac{1}{\Gamma(\theta)} \int_0^1 \frac{1-y}{\alpha B(y)}\bigg(\ln\frac{1}{y}\bigg)^{\theta-1}dy\\
\ar\leq\ar
\frac{1}{\Gamma(\theta)} \int_0^1 \frac{A(y)}{\alpha B(y)}\bigg(\ln\frac{1}{y}\bigg)^{\theta-1}dy<\infty.
 \eeqnn
By (\ref{4.5}) we have
 \beqnn
\psi_i(0,s)=\sum_{k=1}^\infty \bigg(\int_0^\infty p_{ik}(t)dt\bigg) s^k <\infty
 \eeqnn
for each $0\leq s< 1$. Letting $\lambda\to0$ in (\ref{4.17}), we have
 \beqnn
\psi_i(0,s)\leq C_i
+ \frac{1}{\Gamma(\theta)} \int_0^s \frac{A(y)\psi_i(0,y)}{\alpha B(y)}\bigg(\ln \frac{1}{y}\bigg)^{\theta-1}dy.
 \eeqnn
Using the Gronwall's inequality, we have
 \beqlb\label{4.13}
\psi_i(0,s)\leq C_i \exp\bigg[\frac{1}{\Gamma(\theta)}\int^s_0 \frac{A(y)}{\alpha B(y)} \bigg(\ln\frac{1}{y}\bigg)^{\theta-1}dy\bigg].
 \eeqlb
Letting $s\uparrow 1$ we see
 \beqnn
\lim_{s\uparrow 1}\psi_i(0,s) = \lim_{s\uparrow 1}\sum_{j=1}^\infty \int^\infty_0 \tilde{p}_{ij}(t)s^jdt
 =
\int^\infty_0 (1-\tilde{p}_{i0}(t)) dt = E_i (\tilde{\sigma}_0).
 \eeqnn
Hence (\ref{4.15}) follows from (\ref{4.2}) and (\ref{4.13}).

Similarly, by~(\ref{4.12}) we have $$\psi_i(\lambda,s)\geq \frac{1}{\Gamma(\theta)}\int^s_0 \frac{\lambda\tilde{\phi}_{i0}(\lambda)-y^i}{\alpha B(y)}\big(\ln\frac{s}{y}\big)^{\theta-1}dy. $$Letting $\lambda\rightarrow 0$ and then letting $s\rightarrow 1$, we obtain~(\ref{4.16}).\qed

\section{Ergodicity and strong ergodicity}

\setcounter{equation}{0}

One of the main steps to prove Theorems~\ref{t4.1} and~\ref{t5.1} is to compare our nonlinear branching process with immigration with a suitably designed birth-death process, which we now introduce. A similar birth-death process was used by \cite{ChenRR97} in her study of the regularity of the nonlinear branching process with resurrection. Let
 \beqnn
L = M+b_0-1 = \sum_{k=1}^\infty k b_{k+1}
 \eeqnn
and let $(\hat{X}_t)$ be a birth-death process with birth rate $d_i = r_iL+ \gamma m$ and death rate $c_i = r_ib_0$. We denote the density matrix of $(\hat{X}_t)$ by $(\hat{q}_{ij})$. Let $T_0:= \inf\{t\geq 0: \hat{X}_t= 0\}$.

\blemma\label{l5.2} {\rm(1)} Suppose that $m<\infty$, $M< 1$, $r_i$ is increasing and $\sum_{i=1}^\infty r_i^{-1} <\infty$. Then the birth-death process $(\hat{X}_t)$ is strongly ergodic.

{\rm(2)} Suppose that $m<\infty$, $M\leq 1$, $r_i$ is non-decreasing and $\sum_{i=1}^\infty r_i^{-1} <\infty$. Then the birth-death process $(\hat{X}_t)$ is ergodic. \elemma

\proof (1) It is easy to check that the birth-death process is regular. Fix an $\varepsilon> 0$ satisfying $L + \varepsilon< b_0$. Then there exists an $N$ such that $d_i \leq \gamma _i(L+\varepsilon)$ for each $i> N$. Let
 \beqlb\label{5.1+}
S = \sum_{n=1}^\infty\bigg(\frac{1}{c_{n+1}}+ \sum_{k=1}^n \frac{d_k\cdots d_n}{c_k\cdots c_{n+1}}\bigg).
 \eeqlb
It is obvious that $\sum_{n=1}^\infty c_{n+1}^{-1}< \infty$. Notice that for each $n> N$ we have
 \beqnn
\sum_{k=1}^n \frac{d_k\cdots d_n}{c_k\cdots c_{n+1}}\ar\leq\ar \max_{1\leq k\leq N}\frac{d_k\cdots d_N}{c_k\cdots c_N}
\cdot N \cdot \frac{d_{N+1}\cdots d_n}{c_{N+1}\cdots c_{n+1}}+ \sum_{k=1}^{n-N}\frac{d_{N+k}\cdots d_n}{c_{N+k}\cdots c_{n+1}}\cr
\ar\leq\ar N \rho^{n-N}\max_{1\leq k\leq N}\frac{d_k\cdots d_N}{c_k\cdots c_N} + \sum_{k= 1}^{n-N}\frac{\rho^{n-N-k+1}}{c_{n+1}},
 \eeqnn
where $\rho = b_0^{-1}(L +\varepsilon)<1$. Then $S< \infty$. By Corollary 2.4 of \cite{Zhang01},
we conclude that $(\hat{X}_t)$ is strongly ergodic.

(2) Since $L \leq b_0$, we have
 \beqnn
R := \sum_{n=1}^\infty \frac{d_0\cdots d_{n-1}}{c_1\cdots c_n}
 \leq
d_0\sum_{n=1}^\infty\frac{(c_1+\gamma m)(c_2+\gamma m)\cdots (c_{n-1}+\gamma m)}{c_1c_2\cdots c_n}.
 \eeqnn
Taking logarithm on the right-hand side we get
 $$
\ln \bigg(\frac{(c_1+\gamma m)(c_2+\gamma m)\cdots (c_{n-1}+\gamma m)}{c_1c_2\cdots c_n}\bigg)
 =
\sum_{i=1}^{n-1}\ln \bigg(1+\frac{\gamma m}{r_i b_0}\bigg) + \ln \bigg(\frac{1}{c_n}\bigg).
 $$
Since $\lim_{i\rightarrow \infty} r_i=\infty$, we have $\ln{\big(1+\frac{\gamma m}{r_i b_0}\big)} \sim \frac{\gamma m}{r_i b_0}$ as $i\rightarrow\infty$. Then there exists a constant $C\ge 0$ such that for sufficiently large $n$,
 $$
\ln \bigg(\frac{(c_1+\gamma m)(c_2+\gamma m)\cdots (c_{n-1}+\gamma m)}{c_1c_2\cdots c_n}\bigg)
 \leq
C\sum_{i=1}^{\infty}\frac{1}{r_i} + \ln \bigg(\frac{1}{c_n}\bigg),
 $$
and hence
 $$
\frac{(c_1+\gamma m)(c_2+\gamma m)\cdots (c_{n-1}+\gamma m)}{c_1c_2\cdots c_n}\leq \frac{T}{c_n}
 $$
for another constant $T\ge 0$. That implies $R<\infty$. By Theorem 4.55 in \cite{ChenMF04} the birth-death process is ergodic.
\qed

\blemma\label{l4.2} If the nonlinear branching process with immigration has a stationary distribution $\mu=(\mu_j)$, then the generating function $f(s):=\sum_{j=0}^\infty \mu_j s^j$ satisfies the following equation:
 \beqlb\label{n4.1}
\Gamma(\theta)f(s) = \Gamma(\theta)\mu_0+\int_0^s \frac{A(y)}{\alpha B(y)} \bigg(\ln\frac{s}{y}\bigg)^{\theta-1}f(y)dy,\qquad s\in[0,1]. \eeqlb
\elemma

\proof The stationary distribution $(\mu_j)$ satisfies $\mu Q = 0$. In
view of (\ref{1.2}), we have
 \beqlb\label{4.1+}
\quad\mu_j(\gamma + \alpha j^\theta)=\sum_{i=0}^{j-1} \mu_i\gamma
a_{j-i}+\sum_{i=0}^{j+1} \mu_i\alpha i^\theta b_{j-i+1}.
 \eeqlb
Multiplying $s^j$ on both sides of the above equality and then summing over $j$, we have
 \beqlb\label{4.7}\gamma
\sum_{j=1}^\infty \mu_j s^j + \alpha s \sum_{j=1}^\infty \mu_j j^\theta
s^{j-1}
 \ar=\ar
\gamma\sum_{j=1}^\infty \sum_{i=0}^{j-1} \mu_i a _{j-i} s^j
+\alpha\sum^\infty _{j=1} \sum_{i=1}^{j+1} \mu_i i^\theta
b_{j-i+1} s^j. \eeqlb
Interchanging the order of summation,
 \beqlb
\mbox{l.h.s. of } (\ref{4.7})\ar=\ar \gamma \sum_{i=0}^\infty \mu_i s^i \sum_{j=i+1}^\infty a_{j-i} s^{j-i}-
\alpha \mu_1b_0\cr
 \ar\ar
+\,\alpha\sum_{i=1}^\infty \mu_i i^\theta s^{i-1} \sum_{j=i-1}^\infty
b_{j-i+1} s^{j-i+1}\cr
 \ar=\ar
\gamma \sum_{i=0}^\infty \mu_i s^i F(s)- \alpha \mu_1b_0+\alpha
\sum_{i=1}^\infty \mu_i i^\theta s^{i-1}G(s).
 \eeqlb
Letting $j=0$ in (\ref{4.1+}), we see $\mu_0\gamma=\alpha\mu_1b_0$.
Therefore, from (\ref{4.7}) it follows that
 \beqnn \sum_{j=1}^\infty \mu_j j^\theta
s^{j-1} = \frac{f(s)A(s)}{\alpha B(s)}.
 \eeqnn
Multiplying the above equation by $(\ln\frac{s}{y})^{\theta-1}$ and integrating the both sides, we have
 \beqlb\label{3.3}
\int_0^s \sum_{j=1}^\infty \mu_j j^\theta y^{j-1} \bigg(\ln\frac{s}{y}\bigg)^{\theta-1} dy = \int_0^s \frac{A(y)}{\alpha B(y)} \bigg(\ln\frac{s}{y}\bigg)^{\theta-1}f(y) dy.
 \eeqlb
Letting $y=se^{-\frac{x}{j}}$ we get
 \beqnn
\mbox{l.h.s. of } (\ref{3.3})
 \ar=\ar
- \int^\infty _0 \sum_{j=1}^\infty \mu_j j^\theta (s e^{-\frac{x}{j}})^{j-1}\bigg(\frac{x}{j}\bigg)^{\theta-1}
\bigg(-\frac{s}{j}e^{-\frac{x}{j}}\bigg) dx\cr
 \ar=\ar
 \int_0^\infty \sum_{j=1}^\infty \mu_j s^j x^{\theta -1} e^{-x} dx = \Gamma(\theta) [f(s)-\mu_0].
 \eeqnn
Then $f(s)$ is a solution to the differential equation (\ref{n4.1}).\qed

\noindent\emph{Proof of Theorem~{\rm\ref{t4.1}}.}~
(1) By Lemma~\ref{l5.2} the birth-death process $(\hat{X}_t)$ is ergodic. Thus by Theorem 4.45 in \cite{ChenMF04}, the equation
 \beqlb\label{4.11}
u_0=0, ~ d_i(u_{i+1}-u_i)+c_i(u_{i-1}-u_i)+1=0, \qquad i\neq 0
 \eeqlb
has a finite nonnegative solution $(u_i)$. By Remark 2.5 of \cite{Zhang01}, we have
 \beqlb\label{5.2}
u_0=0, ~ u_i= \sum_{k=0}^{i-1}\bigg(\frac{1}{c_{k+1}}+\sum_{j=k+1}^\infty\frac{d_{k+1}\cdots d_j}{c_{k+1}\cdots c_{j+1}}\bigg).
 \eeqlb
It is apparent that $u_i\leq u_{i+1}$. Moreover, we have
 \beqnn
u_{i+1}-u_i= \frac{1}{c_{i+1}}+ \sum_{j=i+1}^\infty \frac{d_{i+1}\cdots d_j}{c_{i+1}\cdots c_{j+1}},
 \quad
u_i- u_{i-1}= \frac{1}{c_i}+\sum_{j=i}^\infty \frac{d_i\cdots d_j}{c_i\cdots c_{j+1}}.
 \eeqnn
Since $d_{i+1}/c_{i+1}< d_i/c_i$ and $1/c_{i+1}< 1/c_i$, it is not hard to show that $u_{i+1}- u_i$ is non-increasing in $i\geq 0$. Coming back to the matrix $Q$, for $i\ge 1$,
 \beqlb\label{4.18}
\sum_{j=0}^\infty q_{ij}u_j
 \ar=\ar
\sum_{j=0}^\infty q_{ij}(u_j-u_i)\cr
 \ar=\ar
c_i(u_{i-1}- u_i)+ r_i\sum_{k=1}^\infty b_{k+1} \sum_{l=1}^k (u_{i+l}- u_{i+l-1})\cr
 \ar\ar\qquad\qquad\qquad\qquad
+ \sum_{k=1}^\infty \gamma a_k \sum_{l=1}^k (u_{i+l}- u_{i+l-1})\cr
 \ar\leq\ar
c_i(u_{i-1}-u_i)+ r_i\sum_{k=1}^\infty k b_{k+1} (u_{i+1}- u_{i}) + \gamma \sum_{k=1}^\infty k a_k(u_{i+1}-u_i)\cr
 \ar=\ar
c_i(u_{i-1}-u_i)+ (r_iL + \gamma m)(u_{i+1}-u_i)\cr
 \ar=\ar
c_i(u_{i-1}-u_i)+ d_i(u_{i+1}-u_i) = - 1.
 \eeqlb
and
 \beqlb\label{5.3}
\sum_{j=1}^{\infty} q_{0j}u_j=\sum_{j=1}^{\infty} q_{0j}(u_j-u_1)=\sum_{j=1}^\infty q_{0j} \sum_{i=1}^j(u_i-u_{i-1})< \sum_{j=1}^\infty q_{0j}j u_1\leq \gamma m u_1< \infty.
 \eeqlb
Then $(u_i)$ is a nonnegative bounded solution to the following equation
 \beqnn
\sum_{j=1}^\infty q_{0j}u_j< \infty, ~ \sum_{j=0}^\infty q_{ij}u_j\leq -1, \qquad i\ge 1.
 \eeqnn
By Theorem 4.45 in \cite{ChenMF04} we know the process is positive recurrent.

(2) Suppose that the process is ergodic. Then letting $s=1$ in (\ref{n4.1}) we get
 \beqnn
\infty> \Gamma(\theta) \sum_{j=1}^\infty \mu_j
 \ge
\int_0^1 \frac{\sum_{j=0}^\infty\mu_j y^j A(y)}{\alpha B(y)}\bigg(\ln\frac{1}{y}\bigg)^{\theta-1}dy
 \geq
\int_0^1 \frac{\mu_0 A(y)}{\alpha B(y)}\bigg(\ln\frac{1}{y}\bigg)^{\theta-1}dy.
 \eeqnn
Since $\mu_0 > 0$, we have (\ref{4.8}). Conversely, suppose that (\ref{4.8}) holds. By the strong Markov property, we have
 \beqnn
E_0(\sigma_0)
 =
E_0(\tau_1) + E_0[E_{X_{\tau_1}}(\sigma_0)]
 =
\frac{1}{q_0} + \sum_{i=1}^\infty \frac{q_{0i}}{q_0}E_i(\sigma_0)
 =
\frac{1}{\gamma} + \sum_{i=1}^\infty a_iE_i(\sigma_0).
 \eeqnn
Using (\ref{4.15}) we have
 \beqnn
E_0(\sigma_0)
 \leq
\frac{1}{\gamma}+ \frac{1}{\gamma\Gamma(\theta)}\bigg[\int_0^1 \frac{A(y)}{\alpha B(y)}\bigg(\ln\frac{1}{y}\bigg)^{\theta-1}dy \bigg]\cdot\exp\bigg[\frac{1}{\Gamma(\theta)}\int^1_0 \frac{A(y)}{\alpha B(y)}\bigg(\ln\frac{1}{y}\bigg)^{\theta-1}dy\bigg].
 \eeqnn
By (\ref{4.8}), the right-hand side is finite. Thus the process is ergodic.

(3) By the assumption, there exists $C> 0$ such that $r_i\geq \frac{C}{b_0-\Gamma}i$ for large enough $i$. Therefore
\beqnn \sum_{j=0}^\infty j q_{ij}\ar=\ar \sum_{k=1}^{\infty} (i+k)r_i b_{k+1}  + \sum_{k=1}^{\infty} (i+k)\gamma a_k +(i-1)r_ib_0-(\gamma  +r_i)i\cr
\ar\leq\ar m-r_i (b_0-\Gamma) \leq m-Ci.
\eeqnn
Applying Corollary 4.49 in \cite{ChenMF04}, we know the process is exponentially ergodic.
\qed

\noindent\emph{Proof of Theorem~{\rm\ref{t5.1}}.}~ (1) Using Lemma~\ref{l5.2}, we see the birth-death process $(\hat{X}_t)$ is strongly ergodic. Let $u_i:= E_i(T_0)$ for $i\geq 0$. Applying Theorem 4.44 and  Lemma 4.48 in \cite{ChenMF04}, we find that  $(u_i)$ is a bounded non-negative solution to equation (\ref{4.11}).
By~(\ref{4.18}) and~(\ref{5.3}),  $(u_i)$ is also a non-negative bounded solution to the following equation
 \beqnn
\sum_{j=1}^\infty q_{0j}u_j< \infty, ~ \sum_{j=0}^\infty q_{ij}u_j\leq -1, \qquad i\ge 1.
 \eeqnn
By Theorem 4.45 in \cite{ChenMF04}, we know the process is strongly ergodic.

(2) Suppose that (\ref{5.1}) holds. Then
 $$
\int_0^1 \frac{A(y)}{\alpha B(y)} \bigg(\ln\frac{1}{y}\bigg)^{\theta-1} dy
 \leq
\gamma\int^1_0 \frac{1}{\alpha B(y)}\bigg(\ln\frac{1}{y}\bigg)^{\theta-1}dy<\infty.
 $$
Letting $i\rightarrow\infty$ in (\ref{4.15}), we get
 $$
\sup_i E_i (\sigma_0)\leq \frac{1}{\Gamma(\theta)}\bigg[ \int^1_0
 \frac{1}{\alpha B(y)}\bigg(\ln\frac{1}{y}\bigg)^{\theta-1}dy\bigg]\cdot
 \exp\bigg[\frac{1}{\Gamma(\theta)}\int_0^1 \frac{A(y)}{\alpha B(y)}
 \bigg(\ln\frac{1}{y}\bigg)^{\theta-1} dy\bigg]<\infty.
 $$
Then by Theorem 4.44 in \cite{ChenMF04} the
process is strongly ergodic.

Conversely, suppose that $X_t$ is strongly ergodic. By Theorem 4.44 in \cite{ChenMF04} and~(\ref{4.16}), we know (\ref{5.1}) holds.

(3) By the strong Markov property, for $i\geq 1$ we have $E_i\sigma_0=\sum_{k=1}^i E_k \sigma_{k-1}.$ Notice that
 $$
E_k \sigma_{k-1}\geq E_k [\mbox{time spent at $k$ until the next jump}]=\frac{1}{r_k+\gamma}.
 $$
Thus $E_i\sigma_0 \geq\sum_{k=1}^i (r_k+\gamma)^{-1}$. By the assumption $\sum_{i=1}^\infty r_i^{-1}=\infty$, we have $\sup_iE_i\sigma_0=\infty$. Applying Theorem 4.44 in \cite{ChenMF04}, we know the process is not strongly ergodic.

\qed

 \noindent\textit{Abstract.}The author would like to thank Professors Mu-Fa Chen, Yong-Hua Mao and Yu-Hui Zhang for their advice and encouragement. I am grateful to the two referees for pointing out a number of typos in the first version of the paper.

\end{document}